\documentclass[12pt]{article}
\usepackage[reqno]{amsmath}
\usepackage{amssymb,amscd,latexsym,paralist}
\usepackage{bm}
\usepackage[all]{xy}

\newcommand{\card}{\mathrm{card}\,}

\newcommand{\spec}{\mathrm{spec}\,}
\newcommand{\Aut}{\mathrm{Aut}}

\newcommand{\Isom}{\mathrm{Isom}}

\newcommand{\tA}{\tilde{A}}

\newcommand{\tG}{\tilde{G}}

\newcommand{\Th}{{\mathcal T}}

\newcommand{\emm}{{\mathfrak{m}}}

\newcommand{\fs}{\mathfrak{s}}

\newtheorem{theorem}{Theorem}
\newtheorem{prop}[theorem]{Proposition}
\newtheorem{cor}[theorem]{Corollary}

\begin{document}
\title{A remark on the structure of torsors under an affine group scheme}
\author{Christopher Deninger}
\date{\ }
\maketitle
\section{Introduction} \label{sec1}
Any affine group scheme $G$ over a field $k$ is the inverse limit over a directed poset $I$ of affine group schemes of finite type over $k$. Let $P$ be a torsor under $G$, i.e. a non-empty affine scheme over $k$ with a $G$-action such that the morphism $G \times_k P \to P \times_k P , (g,p) \mapsto (gp, p)$ is an isomorphism.

\begin{theorem} \label{t2}
 Consider an algebraically closed field $k$ and an affine group scheme $G$ over $k$ written as above as an inverse limit $G = \varprojlim G_i$ over $I$ of algebraic groups $G_i$. Assume that one of the following conditions holds:\\
 i) The poset $I$ is (at most) countable;\\
 ii) The cardinality of $I$ is strictly less than the one of $k$.\\
 Then any $G$-torsor $P$ over $k$ is trivial, i.e. $P (k) \neq \emptyset$.
\end{theorem}

The proof of the first part relies on the following proposition. The second part is an application of a Hilbert Nullstellensatz in infinitely many dimensions.

\begin{prop} \label{t1}
 Any torsor $P$ under an affine group scheme $G$ over a field $k$ is the directed inverse limit of affine schemes of finite type over $k$ with faithfully flat transition maps. Here the directed poset over which the limit is taken, can be chosen to have the same cardinality as $I$, a poset for $G$.
\end{prop}

The theorem implies the following result which motivated the present note.

\begin{cor} \label{t3}
 Let $\Th$ be a neutral Tannakian category over an algebraically closed field $k$. Assume that there is a set $J$ of objects which generate $\Th$ as a tensor category satisfying one of the following two conditions:\\
 i) $J$ is countable;\\
 ii) The cardinality of $J$ is strictly less than the one of $k$.\\
 Then any two fibre functors of $\Th$ over $k$ are isomorphic.
\end{cor}
\section{Proofs} \label{sec2}
Consider an affine group scheme $G = \spec A$ over a field $k$ given as an inverse limit $G = \varprojlim \tG_i$ over a directed poset $I$ of affine group schemes $\tG_i$ of finite type over $k$. Such a representation is always possible by \cite{W2} 3.3 Corollary. Writing $\tG_i = \spec \tA_i$, the projection $G \to \tG_i$ corresponds to a homomorphism of commutative Hopf algebras $\tA_i \to A$ over $k$. The image $A_i$ in $A$ of $\tA_i$ is a sub Hopf algebra of $A$ by \cite{MM} Lemma 4.6~ii). The transition maps between the $\tA_i$'s correspond to inclusion maps in $A$ between the $A_i$'s and the natural Hopf-algebra map
\[
 \varinjlim \tA_i \twoheadrightarrow \varinjlim A_i \subset A
\]
is an isomorphism. Hence we have
\[
 \varinjlim \tA_i \xrightarrow{\sim} \varinjlim A_i = A \; .
\]
Setting $G_i = \spec A_i$ we have $G = \varprojlim G_i$ with the same poset $I$ as before and the coordinate rings $A_i$ now being sub Hopf algebras of $A$. For $j \ge i$ in $I$ the inclusion $A_i \subset A_j$ makes $A_j$ a faithfully flat $A_i$-algebra and $G_j \to G_i$ a faithfully flat morphism by \cite{W2} 14.1 Theorem.

Let $P = \spec B$ be a $G$-torsor. The isomorphism $G \times_k P \xrightarrow{\sim} P \times_k P$ corresponds to an isomorphism of $R$-algebras $A \otimes_k R \cong B \otimes_k R$ where we have set $R = B$. As any $k$-algebra, $R$ is faithfully flat. The image $Q_i$ in $B \otimes_k R$ of $A_i \otimes_k R$ is a finitely generated $R$-sub Hopf algebra of $B \otimes_k R$ and for $j \ge i$ the inclusion $Q_i \subset Q_j$ in $B \otimes_k R$ makes $Q_j$ a faithfully flat $Q_i$-algebra. Choose finitely many generators $q^{(\alpha)}$ of the $R$-algebra $Q_i$ and write them in the form $q^{(\alpha)} = \sum_{\nu} b^{(\alpha)}_{\nu} \otimes r^{(\alpha)}_{\nu}$ with $b^{(\alpha)}_{\nu} \in B$ and $r^{(\alpha)}_{\nu} \in R$. Consider the subalgebra $B_i \subset B$ over $k$ generated by the elements $b^{(\alpha)}_{\nu}$. Since the composed maps
\[
B_i \otimes_k R \hookrightarrow Q_i \otimes_k R \xrightarrow{\text{mult}} Q_i \subset B \otimes_k R
\]
is the inclusion $B_i \otimes_k R \subset B \otimes_k R$ it follows that $B_i \otimes_k R \subset Q_i$ and by construction of $B_i$ also $B_i \otimes_k R = Q_i$ in $B \otimes_k R$. The inclusion maps $Q_i \subset Q_j$ do not necessarily descend to the $B_i$'s. To remedy this, let $\Omega$ be the set of finite subsets $\omega \subset I$ having a maximal element. Denote by $i (\omega) \in \omega$ the (uniquely determined) maximal element of $\omega$. Setting $\omega_1 \le \omega_2$ if $\omega_1 \subset \omega_2$ we get a partial ordering of $\Omega$. It is directed since $I$ is directed: For $\omega_1 , \omega_2 \in \Omega$ set $\omega_3 = \omega_1 \cup \omega_2 \cup \{ k \}$ where $k \in I$ is such that $k \ge i (\omega_1) , k \ge i (\omega_2)$. Then $\omega_3 \in \Omega$ and $\omega_1 \le \omega_3 , \; \omega_2 \le \omega_3$. For every $\omega \in \Omega$ let $B_{\omega}$ be the $k$-algebra in $B$ generated by the $k$-algebras $B_i$ for $i \in \omega$. For $i \in \omega$ we have $i \le i (\omega)$ and hence 
$B_i \otimes_k R = Q_i \subset Q_{i (\omega)} = B_{i (\omega)} \otimes_k R$. It follows that $B_{\omega} \otimes_k R = Q_{i (\omega)}$. In particular $B_{\omega}$ is non-zero. For $\omega_1 \le \omega_2$ the inclusion $B_{\omega_1} \subset B_{\omega_2}$ is faithfully flat since the inclusion of $R$-algebras $B_{\omega_1} \otimes_k R = Q_{i (\omega_1)} \subset Q_{i (\omega_2)} = B_{\omega_2} \otimes_k R$ is faithfully flat. Tensoring the inclusion
\[
 \varinjlim B_{\omega} \subset B
\]
with $R$, we obtain an isomorphism
\[
 \Big( \varinjlim_{\omega} B_{\omega} \Big) \otimes_k R = \varinjlim_{\omega} (B_{\omega} \otimes_k R ) = \varinjlim_{\omega} Q_{i (\omega)} = B \otimes_k R \; .
\]
Hence we have $B = \varinjlim B_{\omega}$. This implies proposition \ref{t1}.

With notations as before we set $P_{\omega} = \spec B_{\omega}$. Since $B_{\omega} \neq 0$ we have $P_{\omega} \neq \emptyset$. For $\omega_1 \le \omega_2$ the morphism $P_{\omega_2} \to P_{\omega_1}$, is faithfully flat and in particular surjective. Since $P_{\omega_1}$ and $P_{\omega_2}$ are of finite type over $k$, the morphism induces a surjection between the sets of closed points $|P_{\omega_2}| \to |P_{\omega_1}|$. If $k$ is algebraically closed the closed points are in natural bijection with the $k$-rational points, and hence the map $P_{\omega_2} (k) \to P_{\omega_1} (k)$ is surjective. Note that $P_{\omega} (k) \neq \emptyset$ by the Nullstellensatz since $P_{\omega} \neq \emptyset$ is of finite type over $k$. If $I$ is countable, then $\Omega$ is countable as well and
\[
 P (k) = \varprojlim_{\omega} P_{\omega} (k)
\]
is an inverse limit of non-empty sets where the transition functions are surjective. Inductively we can choose a cofinal sequence $\omega_1 \le \omega_2 \le \ldots$ in $\Omega$ and it follows that
\[
 P (k) = \varprojlim_{\nu} P_{\omega_{\nu}} (k)
\]
is non-empty. This proves the assertion of theorem \ref{t2} under assumption i).

Note that for uncountable directed posets $I$ an inverse limit of non-empty sets $X_i$ with surjective transition maps $\pi_{ji}$ can well be empty, see e.g. \cite{W1}. On the other hand according to \cite{HM} Proposition 2.7 the inverse limit $\varprojlim X_i$ will be non-empty (and the projections surjective) if the $\pi_{ji}$ are closed surjective continuous maps between quasicompact non-empty $T_1$-spaces $X_i$. For the Zariski topologies on the $P_{\omega} (k)$ all conditions are satisfied except that the transition maps need not map closed sets to closed sets. In fact Hochschild and Moore who are dealing with groups $X_i$ use a suitable coset topology for their arguments. It does not seem to be applicable here. I do not know if torsors for general affine group schemes over an algebraically closed field are trivial.

We now prove the assertion of theorem \ref{t2} under condition ii). For finite $I$ the group $G$ is algebraic and the assertion well known. We may therefore assume that $I$ is infinite. With notations as before, we have $A = \varinjlim_i A_i$ and $B = \varinjlim_{\omega} B_{\omega}$ where the cardinalities of the index sets $I$ and $\Omega$ are the same. The algebra $B_{\omega}$ being finitely generated over $k$ it follows that $B$ is a quotient of the polynomial ring $S = k [T_i , i \in I]$ in infinitely many variables $T_i$ indexed by $I$. Since $B \neq 0$ because of $P \neq \emptyset$, there is an ideal $\fs$ in $S$ with $B \cong S / \fs$ and $\fs \neq S$. We now refer to \cite{L} Theorem from which it follows that for $\card I < \card k$ the ideal $\fs$ has an ``algebraic zero''. This implies that there is a maximal ideal $\emm$ in $S$ with $\fs \subset \emm$ and $S / \emm = k$. The composition
\[
 B \cong S / \fs \longrightarrow S / \emm = k
\]
provides a $k$-rational point of $P$ as claimed.

Generalizations of the Hilbert Nullstellensatz to polynomial rings in infinitely many variables go back to Krull \cite{K} \S\,3, Satz 4 who treats the case of countable $I$ and remarks that the case of arbitrary $I$ could be done along the same lines.

Corollary \ref{t3} is a consequence of theorem \ref{t2} because for two fibre functors $\eta_1$ and $\eta_2$ of $\Th$ the affine $k$-scheme $P = \Isom^{\otimes} (\eta_1 , \eta_2)$ is a torsor under the affine group scheme $G = \Aut^{\otimes} (\eta_1)$. Arguing as in the proof of \cite{DM} Proposition 2.8 one sees that $G$ satisfies the condition of theorem \ref{t2} by our assumption on $\Th$. Hence $P (k) \neq \emptyset$ and thus there exists an isomorphism between $\eta_1$ and $\eta_2$ over $k$. 

\end{document}